\documentclass[a4papeer,20pt]{article}
\usepackage{amsmath, amsthm,graphicx,amsfonts,amssymb,mathrsfs,showidx, multicol}

\usepackage[]{todonotes}

\usepackage{tcolorbox}
\definecolor{mycolor}{rgb}{0.122, 0.435, 0.698}

\newtcbox{\mybox}{on line,
  colframe=mycolor,colback=mycolor!10!white,
  boxrule=0.5pt,arc=4pt,boxsep=0pt,left=6pt,right=6pt,top=6pt,bottom=6pt}

\usepackage{relsize}
\usepackage{color}
\definecolor{rojo}{rgb}{1,0,0}
\usepackage[all]{xy}
\makeindex
\usepackage{hyperref}
\usepackage{color}
\definecolor{abelian}{cmyk}{0.50,0,1,.4}
\definecolor{noabelian}{cmyk}{0.94,0.54,0,0}
\definecolor{rojo}{cmyk}{0,1,1,0}
\definecolor{verde}{cmyk}{0.91,0,0.88,0.12}
\xyoption{arc}

\usepackage{tikz}
\usetikzlibrary{tqft}

\newtheorem{thm}{Theorem}
\newtheorem{prop}[thm]{Proposition}

\theoremstyle{definition}

\theoremstyle{remark}

\newcommand{\op}{\operatorname}

\newcommand{\ben}{\begin{equation}}
\newcommand{\een}{\end{equation}}
\newcommand{\bena}{\begin{equation*}}
\newcommand{\eena}{\end{equation*}}

\newcommand{\Tot}{\longmapsto}

\def\ZZ{\mathbb{Z}}

\def\CC{\mathbb{C}}

\def\RR{\mathbb{R}}

\begin{document}

\title{Thom's counterexamples for the Steenrod problem}
\author{Andres Angel\thanks{Universidad de los Andes, Bogot\'a, Colombia\newline ja.angel908@uniandes.edu.co}, Carlos Segovia\thanks{Instituto de Matem\'aticas, UNAM-Oaxaca, M\'exico\newline csegovia@matem.unam.mx} 
and Arley Fernando Torres\thanks{Universidad Externado de Colombia, Bogot\'a, Colombia\newline arley.torres@uexternado.edu.co}}

\maketitle

\begin{abstract} 
The present paper deals with integral classes $\xi_p\in H_{2p+1}(L^{2p+1}\times L^{2p+1})$ which are counterexamples for the Steenrod realization problem, where $L^{2p+1}$ is the $(2p+1)$-dimensional lens space and $p\geq 3$ is a prime number. For $p=3$, this is Thom's famous counterexample. 
We give a geometric description of this class using the theory of stratifolds.
As a consequence, we obtain a geometric interpretation of the obstruction to realizability in terms of the Atiyah--Hirzebruch spectral sequence.
\end{abstract}

\section*{Introduction}
The Steenrod problem, posted as Problem 25 in \cite{SamEil}, states the following:
``If $z_n\in H_n(K)$ is an $n$-dimensional (integral) homology class in a simplicial complex $K$, does there exist an oriented manifold $M$ and a map $f: M\rightarrow  K$ such that $z_n$ is the image of the generator of $H_n(M)$?\,''

Thom \cite{Th} constructed a counterexample with an integral class $\xi\in H_*(L^{7}\times L^7)$, where $L^7$ is the $7$-dimensional lens space.
This construction extends for $n=2p+1$, with $p$ an odd prime number, by taking the product of two $(2p+1)$-dimensional lens spaces. 
More precisely, the cohomology algebra ($\ZZ_p$-coefficients) has the form $H^*(L^{2p+1};\ZZ_p)\cong\ZZ_p[u,\nu]/(\nu^2=0,u^{p+1}=0)$, with $|\nu|=1$, $|u|=2$ and $u=\beta(\nu)$, where $\beta$ is the Bockstein map. Consider the class
\ben X_p:=u_1\nu_2u_2^{p-1}-\nu_1u_2^p\in H^{2p+1}(L^{2p+1}\times L^{2p+1};\ZZ_p)\,,\een
which we call Thom's counterexamples. The class $X_p$ is an integral class because of the equality $\beta(\nu_1\nu_2u_2^{p-1})=u_1\nu_2u_2^{p-1}-\nu_1u_2^p$.
Now, take the Poincar\'e dual   
\ben
\xi_p\in H_{2p+1}(L^{2p+1}\times L^{2p+1})\,,
\een
and because of the following calculation
\ben    \beta P^1(X_p)=\beta\left(P^1(u_1)\nu_2u_2^{p-1}+u_1\nu_2P^1(u_2^{p-1})-\nu_1P^1(u_2^p)\right)= \beta(u_1^p\nu_2u_2^{p-1})= u_1^pu_2^p\neq 0\,,\een
the integral class $\xi_p$ is a counterexample for the Steenrod realization problem. Notice that the example of Thom is $\xi=\xi_3$.
 
Thom \cite{Th} portrayed the obstruction to realizability regarding cohomology. 
For a solution in terms of homology -as opposed to cohomology- we need a geometric approach for representing the cycles, which is robust enough to treat singularities. We choose the theory of `stratifolds' developed by Kreck in \cite{Kre} (where the stratifolds are compact without boundary and oriented). 
Indeed, we will see that Thom's counterexample can be represented by a $(2p+1)$-dimensional stratifold, where the singular part is a torus with an open neighborhood, which is (homeomorphic to) a cone on $\CC P^{p-1}$.

Conner and Floyd \cite{CF} rephrased the Steenrod problem in terms of the Atiyah--Hirzebruch spectral sequence for oriented bordism $(E_{s,t}^r,d^r_{s,t})$. 
More precisely, the homomorphism from oriented bordism to integral homology $\Omega_*(X)\rightarrow H_*(X)$ is an epimorphism if and only if the differentials 
$d^r_{s,t}:E_{s,t}^r\longrightarrow E_{s-r,t+r-1}^r$ are trivial for all $r\geq 2$. 

Sullivan \cite{Sul} formulated the Steenrod problem to resolve singularities of embedded geometric cycles. 
Indeed, this consists of a blowing-up process to simplify the singularities of geometric cycles, which are constructed inductively by attaching and dragging generalized handles. 
He shows that for a geometric cycle $V$ with singularity $S(V)$, the obstruction for resolving the singularity lies in $H_s(S(V);\Omega_r)$, with $r+s+1=\op{dim}V$.

Stratifolds were introduced by Kreck \cite{Kre} as a generalization of manifolds to provide a positive solution for the Steenrod realization problem.  In the category of CW complexes, there is a stratifold bordism theory $SH_*$ with a natural isomorphism to integral homology. A stratifold $S$ of dimension $n$ is a topological space with a sheaf of functions, which divides the topological space into strata $S_0,\cdots, S_n$, where each stratum $S_i$ is an $i$-dimensional manifold. The singular part of the stratifold $S$ is given by the union $\op{sing}(S):=\bigcup^{n-1}_{i=0}{S_i}$.

There are generators $\alpha_i\in H_i(B\ZZ_p;\ZZ_p)$ such that for $i$ even, $\beta(\alpha_i)$ is the $(i-1)$-dimensional
lens space and for $i$ odd, $\alpha_i$ is the $\operatorname{mod} p$ image of the $i$-dimensional lens space. 
In this respect, Thom's counterexample corresponds to   
\ben
\xi_p:=\alpha_{2}\times \alpha_{2p-1}+\alpha_{1}\times \alpha_{2p} \in  H_{2p+1}(B(\ZZ_p\times \ZZ_p))\,,
\een
where we take the product of $\ZZ_p$-actions to produce and action of $\ZZ_p\times \ZZ_p$. 
The description of the generators is as follows: 
\begin{itemize}
    \item[i)] $\alpha_2$ is a closed surface of genus $(p-1)(p-2)/2$ with an action of $\ZZ_p$ with exactly $p$ fixed points which are removed;
    \item[ii)] $\alpha_n$, for $n$ odd, is represented by $S^n$ with the standard action of $\ZZ_p$ obtaining the $n$-dimensional lens space; and 
    \item[iii)]  $\alpha_{2p}$ is obtained by taking the stratifold, which is the product of $S^1$ with the cone of $\CC P^{p-1}$ equipped with the diagonal action of $\ZZ_p$ which is free. For $\lambda=\op{exp}(2\pi i/p)$, the action of $\ZZ_p$ over $\CC P^{p-1}$ is given by $T([z_1:\dots:z_p])=[z_1:\lambda z_2:\dots ,\lambda^{p-1}z_p]$. 
\end{itemize}

The obstruction to realizability of Thom's counterexample is encoded by the $2p-1$ differential 
\ben
d^{2p-1}:H_{2p+1}(B(\ZZ_p\times \ZZ_p);\Omega_0)\longrightarrow H_{2}(B(\ZZ_p\times \ZZ_p);\Omega_{2p-2})\,,
\een
and the geometric description of the Atiyah--Hirzebruch spectral sequence \cite{Hag}, implies the image $d^{2p-2}(\xi_p)$ is 
\ben\label{678}
[S^1\times S^1\times \CC P^{p-1}\rightarrow X^{2p}]\,,
\een
for $X^{2p}$ the $2p$-skeleton with $X=B(\ZZ_p\times \ZZ_p)$. This element coincides with $M^{2p-2}\times S^1\times S^1$ (where $M^{2p-2}$ is the first Milnor generator of $\Omega/p\Omega$) which is non trivial in $H_2(X;\Omega_{2p-2})$. This shows that $\xi_p$ is not realizable and rectifies the calculations of Koshikawa \cite{koshi}, which were unnoticed by the mathematical community. 

This article is organized as follows: in Section \ref{sec1}, we give a brief review of the Steenrod realization problem.
In Section \ref{georea}, we introduce a basis for the integral homology of $B(\ZZ_p\times\ZZ_p)$, and we show the equivalence of Thom's counterexamples in homology and cohomology.  
Finally, Section \ref{sec2} shows that the image of Thom's counterexample under the $2p-1$ differential of the Atiyah--Hirzebruch spectral sequence is not trivial.

{\bf Acknowledgements:} We thank Instituto de Matem\'aticas-UNAM, Unidad Oaxaca, and Universidad de los Andes for the hospitality and financial support that made this collaboration possible. The first author acknowledges and thanks the hospitality and financial support provided by the Max Planck Institute for Mathematics in Bonn. This work was partially supported by the grant (\#INV-2023-162-2835) from the Fondo de Investigaciones de la Facultad de Ciencias de la Universidad de los Andes. The second author is supported by c\'atedras CONACYT and Proyecto CONACYT ciencias b\'asicas 2016, No. 284621.
The third author would like to thank the Universidad Pontificia Javeriana and especially the Universidad Externado de Colombia, where he has been a professor in the mathematics department since 2020.
\section{The Steenrod problem}\label{sec1}
For $K$ a polyhedron of finite dimension $m$ and $z\in H_*(K)$, we say that 
the integer class $z$ is realizable if there exists a manifold $W$ and a map $f:W\rightarrow K$
such that $z$ is the image of the fundamental class, i.e., $z=f_*([W])$.

Push the class $z\in H_{n-k}(K)$ forward to a smooth manifold $V^n$ along an inclusion $K\hookrightarrow V^n$, constructed by embedding $K$ into the Euclidean space $\RR^n$, with $n\geq 2m+1$.
If $z\in H_{n-k}(V^n)$ is realizable by $f:W^{n-k}\hookrightarrow V^n$, then the Poincar\'e dual $u\in H^k(V^n)$ is the pullback of the universal Thom class in $H^k(MSO(k))$.
The main consequence of this construction is the following.

\begin{thm}[\cite{Th}]\label{Thom}
A necessary condition for an integer cohomology class $x$ to be realizable is that all Steenrod $p$-powers 
$\beta P^i(x)$ 
vanish for all odd prime $p$.
\end{thm}

Novikov gives a specific homology version of the result of Thom as follows.

\begin{thm}[\cite{Nov}]\label{Nov}
Suppose that the integral homology group $H_{n-2i(p-1)-1}(X)$ has no $p$-torsion for every odd prime $p$ and every $i \geq 1$. Then any homology class $z \in H_n(X)$ is realizable.
\end{thm}

Consider the homomorphism $\mu:\Omega_n(X)\rightarrow H_n(X)$, which sends every manifold to the image of its fundamental class.
Conner-Floyd rephrased the Steenrod realization problem in terms of the Atiyah--Hirzebruch spectral sequence $(E_{s,t}^r,d^r_{s,t})$, with the following result.

\begin{thm}[\cite{CF}]\label{CF}
If $X$ is a CW complex, then the differentials $d^r_{s,t}:E^r_{s,t}\rightarrow E^r_{s-r,t+r-1}$
for the Atiyah--Hirzebruch spectral sequence $(E_{s,t}^r, d^r_{s,t})$ are trivial for all $r\geq 2$ if and only if the map  
$\mu:\Omega_n(X)\longrightarrow H_n(X)$ is an epimorphism for all $n\geq 0$.
\end{thm}

To find counterexamples for the Steenrod problem, we notice that the Atiyah--Hirzebruch spectral sequence is trivial modulo the odd torsion part. Hence, every element of $2$-torsion is realizable. Moreover, for a CW complex, if the term  $E^r_{s,t}$ consists entirely of elements of order 2 for $t\neq 0 \op{mod} 4$, then the differential 
$d^r:E^r_{s,t}\rightarrow E^r_{s-r,t+r-1}$ is trivial unless $t=0\op{mod}4$ and $r=1\op{mod} 4$. We will show that also the differentials $d^{4k+1}$ are trivial for $4k+1<2p-2$. The obstruction to realizability is in degree $n-2i(p-1)-1$, for $n=2p+1$, $i=1$, and $p$ an odd prime, hence in dimension $2$. 

\section{Geometric cycles for $B(\ZZ_p\times \ZZ_p)$}
\label{georea}
The Bockstein exact sequence of $B\ZZ_p$ implies the homomorphisms 
$$H_{2n-1}(B\ZZ_p)\stackrel{\op{mod}p}{\longrightarrow} H_{2n-1}(B\ZZ_p;\ZZ_p)\textrm{ and }H_{2n}(B\ZZ_p;\ZZ_p)\stackrel{\beta}{\longrightarrow} H_{2n-1}(B\ZZ_p)$$ are isomorphisms for $n>0$.
Take generators $\alpha_i\in H_i(B\ZZ_p;\ZZ_p)$, such that $\beta(\alpha_i)=\alpha_{i-1}$ for $i$ even, and $\beta(\alpha_i)=0$ for $i$ odd. 
Let $X:=B(\ZZ_p\times \ZZ_p)$. We consider the following commutative diagram between the two Bockstein sequences, 
\ben\label{rel}
\xymatrix{\cdots\ar[r]^(0.4){\times p}&H_{n}(X)\ar[r]^{\op{mod}p}\ar[d]&H_{n}(X;\ZZ_p)\ar[r]^\beta\ar[d]^{=} & H_{n-1}(X)\ar[d]^{\operatorname{mod}p}\ar[r]^(0.6){\times p}&\cdots\\
\cdots\ar[r]^(0.4){\times p}&H_{n}(X;\ZZ_{p^2})\ar[r]&H_{n}(X;\ZZ_p)\ar[r]^{\tilde{\beta}} & H_{n-1}(X;\ZZ_p)\ar[r]^(0.6){\times p}&\cdots\,.}
\een
For $n>1$, the map $\operatorname{mod} p$ is injective, hence every element in the kernel of $\tilde{\beta}$ can be extended to an integral class. 
Thus $H_{2n}(X)$ is generated by $\alpha_{2i-1}\times \alpha_{2n-2i+1}$ ($i=1,\cdots, n$), $H_{2n-1}(X)$ is generated by $\alpha_{2i}\times \alpha_{2n-2i-1}+\alpha_{2i-1}\times \alpha_{2n-2i}$ ($i=0,1,\cdots,n$), and $H_0(X)$ is generated by $\alpha_0\times \alpha_0$.

We represent the integral homology cycles in $H_*(B\ZZ_p)$ by free actions over stratifolds $S$ which is equivalent to singular maps from the quotient $S/\ZZ_p\rightarrow B\ZZ_p$.
The description of the homology cycles $\alpha_i$ is as follows:
\begin{itemize}
    \item[(i)] 
    $\alpha_2$ is the closed oriented surface of genus $(p-1)(p-2)/2$, where there is an action of  $\ZZ_p$ with exactly $p$ fixed points, where we remove around each fixed point a small disc. For example, for $p=3$, this surface is a torus 
    constructed through the lattice $\RR^2$ generated by $e_1$ and the rotation of $e_1$ by 120 degrees. There is an action of $\ZZ_3$ with exactly three fixed points; see Figure \ref{fig1}.
    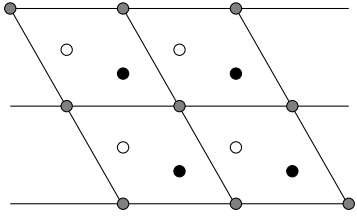
\begin{figure}[h!]
\centering
\begin{tikzpicture}[scale=1.5]

\draw (-1.5,0.866) -- (1.5,0.866);
\draw (-1.5,0) -- (1.5,0);
\draw (-1.5,-0.866) -- (1.5,-0.866);

\draw (-1.5,0.866) -- (-0.5,-0.886);
\draw (-0.5,0.866) -- (0.5,-0.886);
\draw (0.5,0.866) -- (1.5,-0.886);

\draw [fill] (0.5,0.2886) circle [radius=0.05];
\draw  (0,0.5) circle [radius=0.05];

\draw  (-1,0.5) circle [radius=0.05];
\draw [fill] (-0.5,0.2886) circle [radius=0.05];

\draw [fill] (1,-0.5774) circle [radius=0.05];
\draw [fill] (0,-0.5774) circle [radius=0.05];

\draw  (-0.5,-0.366) circle [radius=0.05];
\draw  (0.5,-0.366) circle [radius=0.05];

\draw [fill=gray] (0,0) circle [radius=0.05];
\draw [fill=gray] (1,0) circle [radius=0.05];
\draw [fill=gray] (-1,0) circle [radius=0.05];

\draw [fill=gray] (-0.5,0.866) circle [radius=0.05];
\draw [fill=gray] (0.5,0.866) circle [radius=0.05];
\draw [fill=gray] (-1.5,0.866) circle [radius=0.05];
\draw [fill=gray] (0.5,-0.866) circle [radius=0.05];
\draw [fill=gray] (-0.5,-0.866) circle [radius=0.05];
\draw [fill=gray] (1.5,-0.866) circle [radius=0.05];

\end{tikzpicture}
\caption{The lattice of the torus with three fixed points.}
\label{fig1}
\end{figure}

    \item[(ii)] $\alpha_{2n-1}$ is represented by the sphere $S^{2n-1}=\{(z_1,\cdots,z_n)\in\CC^n:\sum_{i=1}^n|z_i|^2=1\}$, with an action of $\ZZ_p$ induced by the diagonal multiplication $T'(z_1,\cdots,z_n)=(\lambda z_1,\cdots,\lambda z_n)$, where $\lambda=\op{exp}(2\pi i/p)$.
    
    \item[(iii)]  
    $\alpha_{2n}$ is determined by the identity $\beta(\alpha_{2n})=\alpha_{2n-1}$. We use the following equation in the $2n-1$ oriented bordism group $\Omega_{2n-1}(B\ZZ_p)$ from Conner-Floyd \cite[p.~144]{CF}   
    \ben\label{cofo}p\alpha_{2n-1}+[M^4]\alpha_{2n-5}+[M^8]\alpha_{2n-9}+\cdots=0\,, 
    \een
    for $n\geq2$,
where the manifolds $M^{4k}$, $k=1,2,\dots$ are constructed inductively in \cite{CF}. 
Therefore, there is a compact oriented manifold $V^{2n}$, with a free action of $\mathbb{Z}_p$, such that 
\ben
\partial V^{2n} = pS^{2n-1}\sqcup (M^4\times  S^{2n-5} )\sqcup (M^8\times S^{2n-9} )\sqcup \cdots 
\een
Denote by $C(M^{4k})$ the cone on $M^{4k}$ and take the result of gluing $V^{2n}$ with $C(M^4)\times S^{2n-5} \sqcup C(M^8)\times S^{2n-9} \sqcup \cdots$. This construction is a stratifold where the boundary is $pS^{2n-1}$. Therefore, the Bockstein of the corresponding homology class is $\alpha_{2n-1}$. Description of integral homology cycles in terms of stratifolds and the Bockstein homomorphism can be found in the book of Kreck \cite{Kre} and the paper by the authors \cite{ATS}. 

For $n=p$ an odd prime number, we can represent $\alpha_{2p}$ using the equation in $\Omega_{2p-1}(B\ZZ_p)$ from Conner-Floyd \cite[p.~95]{CF}
    \ben\label{co3}
      pS^{2p-1}=S^1\times \CC P^{p-1}       
    \een
where the action of $\ZZ_p$ on $S^1$, $S^{2p-1}$ and $\CC P^{p-1}$ are defined for $\lambda=\op{exp}(2\pi i/p)$, by $T_1(z)= \lambda z$, $T'(z_1,\dots,z_{p})=(\lambda z_1,\lambda z_2,\lambda^2 z_3,\dots,\lambda^{p-1}z_p)$ and $T([z_1:\cdots:z_p])= [z_1:\lambda z_2:\cdots :\lambda^{p-1}z_p]$, respectively.
Thus, we consider the product of $S^1$ with the cone of $\CC P^{p-1}$, equipped with the diagonal action of $\ZZ_p$. From \eqref{co3}, the boundary satisfies $\beta(\alpha_{2p})=\alpha_{2p-1}$. 
\end{itemize}

\begin{prop}\label{tryout}
    The generator $\alpha_{2n}$ is bordant to a stratifold with singular part of dimension $2(n-p)+1$ for $n\geq p$ and empty otherwise.
\end{prop}
\begin{proof}
The manifolds $M^{4l}$, with $4l < 2p-2$, belong to $p\Omega_*$; see the paper by Floyd \cite[p. 336]{Floyd} or the book of Conner-Floyd \cite[p. 93]{CF}. Therefore, there exist manifolds $M_l\in \Omega_{4l}$ such that $M^{4l} = pM_l$ for $4l < 2p-2$. 
Thus, the expression \eqref{cofo} has the form 
\begin{equation}\label{Coni2}
p\alpha_{2n-1}+[M_1] p\alpha_{2n-5}+[M_2]p\alpha_{2n-9}+\cdots +[M^{2p-2}]\alpha_{2(n-p)+1}+\cdots=0\,.
\end{equation}
We use again equation \eqref{cofo}, to substitute $p\alpha_{2n-5}$, $p\alpha_{2n-9}$, $\dots$, $p\alpha_{2(n-p)+5}$ in terms of the manifolds $M^4$, $M^8$, $\dots$, $M^{2p-2}$, $\dots$. 
We apply this process inductively by substituting all the elements of the form $p\alpha_{2i-1}$, where if $2i-1< 2p-2$, from Conner-Floyd \cite[p. 95]{CF}, then $\alpha_{2i-1}$ is of order $p$ in $\widetilde\Omega_*(B\ZZ_p)$. Therefore, we can cap these terms with the corresponding bordisms. Thus, the top dimension of the singular part is the term that comes with $M^{2p-2}$, which is $\alpha_{2(n-p)+1}$.
\end{proof}

Thom's counterexample is
\ben
\xi_p:=\alpha_{2}\times \alpha_{2p-1}+\alpha_{1}\times \alpha_{2p} \in  H_{2p+1}(X)\,.
\een
Recall the cohomology of the lens space is 
$H^*(L^{2p+1};\ZZ_p)\cong\ZZ_p[u,\nu]/(\nu^2=0,u^{p+1}=0)$, with $|\nu|=1$, $|u|=2$ and $u=\beta(\nu)$. The Poincar\'e duality correspondence is given by $\alpha_1\longleftrightarrow u^p$ and $\alpha_2\longleftrightarrow \nu u^{p-1}$, which agrees with the Bocksteins
$\beta(\alpha_2)=\alpha_1$ and $\beta(\nu u^{p-1})=u^p$ and conversely, we have $u\longleftrightarrow \alpha_{2p-1}$ and $\nu \longleftrightarrow\alpha_{2p}$, which also agrees with the Bocksteins
$\beta(\alpha_{2p})=\alpha_{2p-1}$ and $\beta(\nu)=u$. Therefore, we apply the Poincar\'e duality isomorphism $D$ to $X_p=u_1\nu_2u_2^{p-1}-\nu_1u_2^p$ as follows:
\begin{align}
D(X_p)&=[L_1^{2p+1}\times L_2^{2p+1}]\cap (u_1\times\nu_2u_2^{p-1})- [L_1^{2p+1}\times L_2^{2p+1}]\cap(\nu_1\times u_2^p)\\
& = D(u_1)\times D(\nu_2u_2^{p-1})+D(\nu_1)\times D(u_2^p)\\
& = \alpha_{2p-1}\times \alpha_2+\alpha_{2p}\times \alpha_1=\xi_p\,.
\end{align}

\section{Stratifolds and the Atiyah--Hirzebruch spectral sequence}
\label{sec2}

Using the notion of stratifolds, many constructions in algebraic topology become simpler. For instance, the Atiyah--Hirzebruch spectral sequence for oriented bordism $(E^r_{s,t},d^r_{s,t})$, has a geometric description \cite{Hag}, using a Postnikov tower $\Omega^{(r)}$. This description of the spectral sequence is similar  to the Conner-Floyd spectral sequence appearing in equivariant bordism \cite{CF} and the spectral sequence for orbifold cobordism of \cite{angel2007spectral}. The bordism theory $\Omega^{(r)}$ is composed of oriented p-stratifolds, with all strata of codimension $0<k<r+2$ empty. 
Thus a stratifold $S$ in $\Omega^{(r)}_n$ is an $n$-dimensional stratifold with singular part of dimension at most $(n-r-2)$. We put a similar restriction on the stratifold bordisms, which are $(n+1)$-dimensional stratifolds with boundary, and the singular part is of dimension at most $(n-r-1)$. 

If we denote by $\Omega_n$ the $n$th coefficient group $\Omega_n(*)$, we have maps $\Omega_n\rightarrow\Omega_n^{(r)}$ which are isomorphisms for $n\leq r$, and $\Omega_n^{(r)}$ is trivial for $n>r$. Among other properties, 
for a CW complex $X$ with $k$-skeleton $X^k$, we obtain that $\Omega_n^{(r)}(X^k)$ is trivial for $k+r< n$. 

For $r\geq 2$, there are natural isomorphisms \ben E^{r}_{s,t}\cong\operatorname{Im}( \Omega^{(t+r-2)}_{s+t}( X^s ) \longrightarrow 
 \Omega^{(t)}_{s+t}( X^{s+r-1} ))\,,\een and the differential $d^r_{s,t}:E^r_{s,t}\longrightarrow E^r_{s-r,t+r-1}$ is induced by the following commutative diagram 
 \ben
\xymatrix{&\Omega_{s+t}^{(t+r-2)}(X^s)\ar[r]\ar[d]_\Phi&\Omega_{s+t}^{(t)}(X^{s+r-1})\ar[d]^\Phi\\
&\Omega_{s+t-1}(X^{s-r+1})\ar[d]\ar[r]&\Omega_{s+t-1}(X^{s-1})\ar[d]\\
\Omega_{s+t-1}^{(t+2r-3)}(X^{s-r})\ar[r]&\Omega_{s+t-1}^{(t+2r-3)}(X^{s-r+1})\ar[r]&\Omega_{s+t-1}^{(t+r-1)}(X^{s-1})
\,,}
\een 
where $\Phi$ is a natural transformation defined by 
\ben
\Omega_n^{(r)}(X)\rightarrow\Omega_n^{(r)}(X,X^{n-r-1})\stackrel{\cong}{\rightarrow}\Omega_n(X,X^{n-r-1})\rightarrow\Omega_{n-1}(X^{n-r-1})\,.
\een
The isomorphism $\Omega_n^{(r)}(X,X^{n-r-1})\stackrel{\cong}{\rightarrow}\Omega_n(X,X^{n-r-1})$ is given by the restriction to the top stratum and the map $\Omega_n(X,X^{n-r-1})\rightarrow\Omega_{n-1}(X^{n-r-1})$ is the boundary homomorphism. Therefore, for a
stratifold $S$ of dimension $s+t$, with a map $f:S\rightarrow X^s$, the image of the differential $d^r_{s,t}$ is induced by
\ben \label{wef}[f:S\rightarrow X^s]\Tot [f_{\op{sing(S)}}\circ g:\partial W\rightarrow X^{s-1}]\,,\een
where $W$ is the top stratum of $S$ and $g:\partial W\rightarrow\op{sing}(S)$ is the attaching map used to glue $W$ to the singular part $\op{sing}(S)$.

\begin{thm}
For $X=B(\ZZ_p\times \ZZ_p)$, the differentials $d^r_{s,t}$ are trivial for $r\leq 2p-2$. In particular, the differential $d_{2p+1,0}^{2p-1}$ is of the form 
$$d_{2p+1,0}^{2p-1}:H_{2p+1}(B(\ZZ_p\times \ZZ_p);\Omega_0)\longrightarrow H_{2}(B(\ZZ_p\times \ZZ_p);\Omega_{2p-2})\,,$$
and the image of the class $\xi_p:=\alpha_{2}\times \alpha_{2p-1}+\alpha_{1}\times \alpha_{2p}$, under the differential $d^{2p-1}$ is non trivial.
\end{thm}

\begin{proof}
        Recall that $H_{2n}(X)$ is generated by $\alpha_{2i-1}\times \alpha_{2n-2i+1}$ ($i=1,\cdots, n$), $H_{2n-1}(X)$ is generated by $\alpha_{2i}\times \alpha_{2n-2i-1}+\alpha_{2i-1}\times \alpha_{2n-2i}$ ($i=0,1,\cdots,n$), and $H_0(B(X)$ is generated by $\alpha_0\times \alpha_0$. 
    We can restrict to differentials $d^r_{2n+1,0}$ since those starting on coordinates $(2n,0)$ are always trivial (the classes $\alpha_{2n+1}$ are always represented by spheres). 
    For $n\geq p$, Proposition \ref{tryout} implies that the singular part of $\alpha_{2i}\times \alpha_{2n-2i-1}+\alpha_{2i-1}\times \alpha_{2n-2i}$ is of dimension $2(n-p)$. The differential $d^r:E^r_{2n-1,0}\rightarrow E^r_{2n-r-1,r-1}$ is of the form 
    \ben
\xymatrix{\operatorname{Im}\left( \Omega^{(r-2)}_{2n-1}\left( X^{2n+1} \right) \longrightarrow 
 \Omega^{(0)}_{2n-1}\left( X^{2n-r-2} \right)
\right)\ar[d]_{d^{2p-1}}\\
\operatorname{Im}\left( \Omega^{(2r-3)}_{2n-2}\left( X^{2n-r-1} \right) \longrightarrow 
 \Omega^{(2n-2)}_{r-1}\left( X^{2n-2} \right)
\right)\,,}
\een
For $r\leq 2p-2$, the element $\alpha_{2i}\times \alpha_{2n-2i-1}+\alpha_{2i-1}\times \alpha_{2n-2i}$ belong to $\Omega^{(r-2)}_{2n-1}\left( X^{2n+1} \right)$ since $2(n-p)\leq 2n-r-1$. The application of $d^r$ to this element is trivial since $2(n-p)\leq 2n-r-2$, and hence we can cone starting from the Milnor generator $M^{2p-2}$. In fact, we have $E^2\cong \dots\cong E^{2p-1}$ since $\Omega_*$ has no odd torsion, we have the commutative diagram 
        \begin{equation}\label{fre}\xymatrix{ E_{s,0}^r\otimes \Omega_t\ar[r] \ar[d]_{d^r\otimes\op{id}}& E^r_{s,t} \ar[d]_{d^r}\\
    E_{s-r,r-1}^r\otimes \Omega_t\ar[r] & E^r_{s-r,t+r-1}\,, }\end{equation}
as in Conner-Floyd \cite[p. 17,41]{CF}, we have by induction that the rows are isomorphisms for $r\leq 2p-2$. 
Finally, for $r=2p-1$, the element $d^{2p-1}_{2p+1,0}(\xi_p)=M^{2p-2}\times S^1\times S^1$ is not zero  in $H_{2}(B(\ZZ_p\times \ZZ_p);\Omega_{2p-2})$, since $M^{2p-2}$ is a Milnor generator of $\Omega/p\Omega$. For $p=3$, 
$M^4$ can be taken to be $\mathbb{CP}^2$ and we find the obstruction to realizability with $d^5$.
\end{proof}

\bibliographystyle{amsalpha}
\bibliography{biblio}
\end{document}